\documentclass[12pt]{article}

\newcommand{\sch}{Schr\"{o}dinger}

\newcommand{\ssz}{\scriptsize}

\newcommand{\np}{\newpage}

\newcommand{\be}{\begin{equation}}
\newcommand{\ee}{\end{equation}}
\newcommand{\bea}{\begin{eqnarray}}
\newcommand{\eea}{\end{eqnarray}}

\newcommand{\hs}{\hspace}
\newcommand{\vs}{\vspace}

\newcommand{\ra}{\rightarrow}

\newcommand{\f}{\frac}     
\newcommand{\nn}{\nonumber}

\newcommand{\al}{\alpha}

\newcommand{\iny}{\infty}

\newcommand{\lab}{\label}
\usepackage{amsmath}
\usepackage{graphicx}
\begin{document}
  
\begin{center}
 {\Large\bf       Bessel Type Orthogonality For Hermite Polynomials  }
\end{center}
\begin{center}
Omid Hamidi 
\end{center}
\begin{center}
{\ssz Department of Physics, Shahid Bahonar University of Kerman, Kerman, IRAN.\\
(e-mail: hamidi@mail.uk.ac.ir)}
\end{center}
\begin{abstract}
It is shown that Hermite polynomials satisfy a Bessel type orthogonality relation based on the zeros of a single indexed Hermite polynomial  and with a finite integration interval. Because of the role of non-symmetric zeros in the final relation, its applicability   covers Hermite polynomials $P_n(x)$ with $n\ge 3$.
\vs{0.5cm}

\end{abstract}
 
\np 

\section{Introduction}
Among the orthogonal functions (and even the polynomials) appearing in mathematical
 physics the Bessel functions $J_{\nu}(x)$ have a different type of orthogonality relation \cite{ak}, 
\be 
\int_0^a dx\; x J_{\nu}(\f{\al_{\nu m} x}{a})J_{\nu}(\f{\al_{\nu n} x}{a})=0\;.
\ee
As seen, the two Bessel functions under the integral have the same function index $\nu$ and the orthogonality is relative to the location of the zeros $\al_{\nu k}$ satisfying $J_{\nu}(\al_{\nu k})=0$. Also the interval, $[0,a]$ ($a$ is arbitrary), over which the integral is carried out is only part of the domain over which $J_{\nu}(x)$ is defined. In this short note it is shown that the Hermite polynomials $H_n(x)$ which satisfy the standard orthogonality relation
\be 
\int_{-\iny}^{\iny}dx\;e^{-x^2}\,H_n(x)H_m(x)=0\;,\hs{5mm}n\ne m
\ee
also satisfy an orthogonality relation similar in spirit to the orthogonality relation of the Bessel functions. In the next section the details of work are presented.
\section{Details Of The Procedure}
To obtain the desired orthogonality relation
consider the Hermite differential equation,
\be 
H_n''(x) -2x H_n'(x)+2n H_n(x)=0\;,\hs{5mm}-\iny <x<\iny
\lab{ehd}
\ee
where as usual prime indicates differentiaion with respect to the argument.
Let us first change the independent variable to $y$ according to $x=k_n y$. 
Equation (\ref{ehd}), upon multiplying by $k_n^2$, then transforms into,
\be 
\f{d^2H_n(k_ny)}{dy^2}- 2k_n^2 y\f{dH_n(k_ny)}{dy}+2nk_n^2 H_n(k_ny)=0
\lab{fss}
\ee
Now in (\ref{fss}), a transformation of the dependent variable   according to 
\be 
H_n(k_ny)=e^{k^2_n y^2/2}  \psi_n(k_ny)
\lab{shm}
\ee
will transform it into
\be 
 \f{d^2 \psi_n(k_ny)}{dy^2}+\left[k_n^2(1+2n)-k_n^4 y^2\right] \psi_n(k_ny)=0
 \lab{dsr}
\ee
where the common exponential factor has been dropped. From (\ref{shm}) it is  noticed that $\psi_n(y)$ is either even or odd like $H_n(y)$ and
  the location of zeros of the $\psi_n(y)$ are symmetric with respect to the origin and on the finite $x$ axis coincide with the zeros of the Hermite polynomail $H_n(y)$. 
 Now let us choose $ k_n$ as follows
\be 
k_n=\f{\al_{nj}}{b}
\ee
where $\al_{nj}$ indicates the location of the $j$th zero of $\psi_n(y)$, that is $\psi_n(\al_{nj})=0$.
The next step is to write (\ref{dsr}) for two different values of  $k_n$, i.e.,
\bea 
 \f{d^2 \psi_n(\f{\al_{ni}y}{b})}{dy^2}&=&\left[\f{\al_{ni}^4 y^2}{b^4}-\f{\al^2_{ni}(1+2n)}{b^2} \right] \psi_n(\f{\al_{ni}y}{b})\nn\\
  \f{d^2 \psi_n(\f{\al_{nj}y}{b})}{dy^2}&=&\left[\f{\al_{nj}^4 y^2}{b^4}-\f{\al^2_{nj}(1+2n)}{b^2} \right] \psi_n(\f{\al_{nj}y}{b})
\eea
Now multiply the first equation by $\psi_n(\f{\al_{nj}y}{b})$ and the second one by 
$\psi_n(\f{\al_{ni}y}{b})$ and subtract the second one from the first and integrate over $y$ from $y=-b$ to $y=b$. Once done, on the Left Side (LS) (Right Side (RS)), one would get
\bea 
\mbox{LS}&=&\int_{-b}^b dy\left[
\psi_n(\f{\al_{nj}y}{b})\f{d^2 \psi_n(\f{\al_{ni}y}{b})}{dy^2}-
\psi_n(\f{\al_{ni}y}{b})\f{d^2 \psi_n(\f{\al_{nj}y}{b})}{dy^2}\right]
\lab{fyd}\\
\mbox{RS}&=&\f{(\al_{ni}^2-\al_{nj}^2)}{b^2} \int_{-b}^b dy
\left[ \f{(\al_{ni}^2+\al_{nj}^2)y^2}{b^2}-(1+2n)\right]\psi_n(\f{\al_{nj}y}{b})\psi_n(\f{\al_{ni}y}{b})
\eea
Due to the formal self-adjointness of the differential operator $\f{d^2}{dy^2}$, the expression LS equals a boundary term given by
\be  
\mbox{LS} =  \left.\left[
\psi_n(\f{\al_{nj}y}{b})\f{d  \psi_n(\f{\al_{ni}y}{b})}{dy}-
\psi_n(\f{\al_{ni}y}{b})\f{d  \psi_n(\f{\al_{nj}y}{b})}{dy}\right]\right|_{y=-b}^{y=b}=0 
\lab{vfd}
\ee
where the fact $\psi_n(\al_{nj})=\psi_n(-\al_{nj} ) =0$ (the same also for $j\ra i$) has been used.
On the right side RS, the coefficeint in front, i.e., $(\al_{ni}^2-\al_{nj}^2)$, will be non-zero if 
$\al_{ni}$ and $\al_{nj}$ are non-symmetric zeros. Therefore, the   result bellow holds   for $n\ge 3$ in which case non-symmetric zeros exist.
So for non-symmetric zeros, $\al_{ni}\ne \al_{nj}$, the equality LS$=$RS will yield
\be 
0=\int_{-b}^b dy
 \left[ \f{(\al_{ni}^2+\al_{nj}^2)y^2}{b^2}-(1+2n)\right]\psi_n(\f{\al_{nj}y}{b})\psi_n(\f{\al_{ni}y}{b})
 \lab{vdg}
\ee
One last step in  (\ref{vdg}) is to use (\ref{shm}) to express $\psi_n(\f{\al_{nj}y}{b})$ in terms of the corresponding Hermite polynomial $H_n(\f{\al_{nj}y}{b})$   and will yield
\be 
0=\int_{-b}^b dy
 \left[ \f{(\al_{ni}^2+\al_{nj}^2)y^2}{b^2}-(1+2n)\right]
e^{-(\al_{ni}^2 +\al_{nj}^2)y^2/2b^2} H_n(\f{\al_{ni}y}{b})H_n(\f{\al_{nj}y}{b})
\lab{gdh}
\ee 
which is the desired orthogonality relation. As it may have been noticed the interval of integration could have been chosen $[0,b]$ in which case LS in  (\ref{vfd}) would be zero again,
\be  
\mbox{LS} =  \left.\left[
\psi_n(\f{\al_{nj}y}{b})\f{d  \psi_n(\f{\al_{ni}y}{b})}{dy}-
\psi_n(\f{\al_{ni}y}{b})\f{d  \psi_n(\f{\al_{nj}y}{b})}{dy}\right]\right|_{y=0}^{y=b}=0 
\lab{hfv}
\ee
where the vanishing of the expression at lower end of the interval  is due to the vanishing of $\psi_n(y)$ or its derivative at $y=0$ (depending on $n$ being odd or even, respectively). So one would get
\be 
0=\int_{0}^b dy
 \left[ \f{(\al_{ni}^2+\al_{nj}^2)y^2}{b^2}-(1+2n)\right]
e^{-(\al_{ni}^2 +\al_{nj}^2)y^2/2b^2} H_n(\f{\al_{ni}y}{b})H_n(\f{\al_{nj}y}{b})
\lab{hfh}
\ee
Having shown the new orthogonality relations (\ref{gdh}) and (\ref{hfh}) it is a question in which, if any, occasion  would such a type of orthogonality play a role. Looking back at (\ref{dsr}) one notices that it is the {\sch} equation for the simple harmonic oscillator problem. In that equation
$k_n^2(1+2n)$ plays the role of the eigenvalue while $k_n^4 y^2$ is the potential energy. Because of the multiplicative factor $k^4_n$ (which also appears in the energy eigenvalue), the potential energy depends on the energy eigenvalue.
Such type of potentials are not trivial in proving the orthogonality of eigenfunction.
  Problems of this sort   have been discussed in the literature \cite{sz,yl,fvg}.


\begin{thebibliography}{99}
\bibitem{ak}G. B. Arfken, H. J. Weber: {\em Mathematical Methods For Physicists},   
5th ed., Academic Press, London  2001.
%
\bibitem{sz}
H. Jallouli, H. Sazdjian, Ann. Phys. {\bf 253} (1997) 376. 
%
\bibitem{yl}R Yekken,     R J Lombard  
J. Phys. A: Math. Theor. {\bf 43} (2010) 125301 (17pp)
%
\bibitem{fvg}
J. García-Martínez, J. García-Ravelo  , J.J. Peña   , A. Schulze-Halberg , 
Physics Letters A {\bf 373} (2009) 3619-3623

\end{thebibliography}
\end{document}